\documentclass{amsart}
\usepackage{amssymb}
\usepackage{color}

\usepackage{a4}
\begin{document}
\newtheorem{theorem}{Theorem}[section]
\newtheorem{lemma}[theorem]{Lemma}
\newtheorem{remark}[theorem]{Remark}
\newtheorem{definition}[theorem]{Definition}
\newtheorem{corollary}[theorem]{Corollary}
\newtheorem{example}[theorem]{Example}
\makeatletter
 \renewcommand{\theequation}{%
  \thesection.\alph{equation}}
 \@addtoreset{equation}{section}
 \makeatother
\def\id{\operatorname{Id}}
\def\qedbox{\hbox{$\rlap{$\sqcap$}\sqcup$}}
\def\ffrac#1#2{{\textstyle\frac{#1}{#2}}}
\def\Tr{\operatorname{Tr}}
\def\nn{\nonumber}
\def\bea{\begin{array}}
\def\eea{\end{array}}
\def\beq{\begin{eqnarray}}
\def\eeq{\end{eqnarray}}
\def\ctw{c_2 (\theta , m)}
\def\cth{c_3 (\theta , m)}
\def\cfo{c_4 (\theta , m)}
\def\cfi{c_5 (\theta , m)}
\def\csi{c_6 (\theta , m)}
\def\cse{c_7 (\theta , m)}
\def\dtw{d_2 (\theta , m)}
\def\dth{d_3 (\theta , m)}
\def\dfo{d_4 (\theta , m)}
\def\don{d_1 (\theta , m)}
\def\gt{\tilde\gamma}
\def\Th{\theta}
\def\pip{\Pi_+}
\def\pim{\Pi_-}
\def\pipl{\Pi_+^\star}
\def\pimi{\Pi_-^\star}
\def\relsum{\sum_{n=0}^\infty d_n (m)}
\newcommand{\nats}{\mbox{${\rm I\!N }$}}
\newcommand{\reals}{\mbox{${\rm I\!R }$}}
\title[Heat Kernel Coefficients for Chiral Bag Boundary Conditions]
{Heat kernel coefficients for chiral bag boundary conditions}
\author{Giampiero Esposito, Peter Gilkey and Klaus Kirsten}
\begin{address}{GE: Istituto Nazionale di Fisica Nucleare,
Sezione di Napoli, Complesso Universitario di Monte S. Angelo, Via
Cintia, Edificio N', 80126 Napoli, Italy}\end{address}
\begin{email}{giampiero.esposito@na.infn.it}\end{email}
\begin{address}{PG: Math. Dept., University of Oregon,
Eugene, Or 97403, USA}\end{address}
\begin{email}{gilkey@darkwing.uoregon.edu}\end{email}
\begin{address}{KK: Department of Mathematics, Baylor University \\
Waco, TX 76798, USA}\end{address}
\begin{email}{Klaus\_Kirsten@baylor.edu}\end{email}
\begin{abstract}
We study the asymptotic expansion of the smeared $L^{2}$-trace
of $f \; e^{-tP^2}$ where $P$ is
an operator of Dirac type, $f$ is an auxiliary smooth
smearing function which is used to localize the problem, and
chiral bag boundary conditions are imposed. Special case calculations,
functorial methods and the theory of $\zeta$- and $\eta$-invariants
are used to obtain the boundary part of the heat-kernel
coefficients $a_{1}$ and $a_{2}$.
\end{abstract}
 \maketitle

\section{Introduction}

Local boundary conditions for operators of Dirac type have been
studied in the physical and mathematical literature with a variety
of motivations over many years. Some key points in this respect
are as follows.
\vskip 0.3cm
\noindent
(i) Local boundary
conditions for massless fermionic fields ruled by a Dirac operator
can be applied to one-loop quantum cosmology \cite{deat91},
\cite{espo92} and are part of the investigation of conformal
anomalies \cite{moss94} in Euclidean field theory \cite{espo94}.
Moreover, they are the first step towards analyzing boundary
counterterms in supergravity theories, with the associated
unsettled issue of proving finiteness \cite{deat96} or lack of
finiteness \cite{dese99} of supergravity theories on manifolds
with boundary. In other words, the local boundary conditions for
fermionic fields are part of a general scheme \cite{luckock91}
leading to locally supersymmetric boundary conditions for
fermionic and bosonic fields \cite{deat91}, \cite{avra99}, and
hence can be used to test perturbative consistency of supergravity
models in cosmological \cite{espo92}, \cite{espo97} or
field-theoretical backgrounds.
\vskip 0.3cm
\noindent
(ii) Local boundary conditions of chiral bag type are a substitute for
introducing small quark masses to drive the breaking of chiral
symmetry \cite{wipf95}. One of the first papers where the chiral
boundary conditions were introduced is the work by Hrasko and Balog
\cite{hras84}, and one of the first applications to chiral bag
models is presented in \cite{gold83}.
\vskip 0.3cm
\noindent
(iii) Chiral bag
boundary conditions have been recently proved to lead to a
strongly elliptic boundary-value problem for the squared Dirac
operator \cite{bene03-36-11533}, and the associated global
heat-kernel asymptotics has been investigated in detail, on the
Euclidean ball, in \cite{espo02-66-085014}. An early paper on the
role of boundary conditions for Dirac operators is in the
framework of fermionic billiards \cite{anto90}, studied even
earlier by Berry and Mondragon \cite{berry87}.
\vskip 0.3cm
For more general Riemannian manifolds with boundary, however, the
investigation of such a global asymptotics in the chiral bag case
is, to our knowledge, an open research field, and it appears
desirable to understand how far can one go by exploiting
functorial methods (e.g. conformal rescalings of the metric) and
special case calculations, which are tools frequently used in
invariance theory \cite{gilk95b}, \cite{kirs01b}. For this
purpose, both algorithms are exploited in our paper, where the
general mathematical setting is as follows.

Let $m=2{\overline m}$ be even and let
$P=\gamma_{j}\nabla_{j}+\psi$ be an operator of Dirac type on a
compact oriented Riemannian manifold $M$ of dimension $m$, where
$\nabla$ is a compatible unitary connection, i.e.
$\nabla\gamma=0$. The spinor space has then dimension
$d_s=2^{{\overline m}}$, and the $\gamma$-matrices can be taken to
be skew-adjoint and obeying the Clifford relation
$$\gamma_i \gamma_j + \gamma_j \gamma_i = -2\delta _{ij}.$$ Near
the boundary, let $e_{m}$ be the inward-pointing unit normal and
$\gamma_{m}$ be the projection of the $\gamma$-matrices on
$e_{m}$. Moreover, the generalization of $\gamma_{5}$ to arbitrary
even dimension is provided by
\begin{equation}
{\widetilde \gamma} \equiv i^{{\overline m}} \gamma_{1}{\ldots}
\gamma_{m}.
\end{equation}
The squared Dirac operator is studied with local boundary
conditions of chiral bag type. These boundary conditions involve a
real angle $\theta$ and they read
\begin{equation}
\left.  \Pi _- \varphi \right|_{\partial M} =0 ,\end{equation}
where we have introduced the `projectors' \beq \Pi_\mp \equiv \frac 1 2
\left( 1\pm  e^{\theta \gt} \gt \gamma_m \right).
\label{projpm}\eeq Under the above assumptions, the squared
operator $P^{2}$ is an operator ${\widetilde P}$ of Laplace type
\cite{gilk95b}. The associated heat kernel can be defined as the
solution, for $t>0$, of the heat equation
\begin{equation}
\left({\partial \over \partial t}+{\widetilde P}\right)
U(x,x';t)=0,
\end{equation}
obeying the initial condition
\begin{equation}
\lim_{t \rightarrow 0}\int_{M}dx' U(x,x';t) \Psi(x')=\Psi(x),
\end{equation}
jointly with the boundary conditions ${\mathcal B}_\theta$ defined
by
\begin{equation}
\left.\Pi_- U(x,x';t)  \right|_{x\in\partial M} =0 , \nn\\
\left. \Pi _- P_x U (x,x'; t)\right|_{x\in\partial M} =0 .
\end{equation}
Here, $dx'$ denotes the Riemannian volume element of the manifold
$M$ and $P_x$ denotes the Dirac operator with respect to the
variable $x$. The $L^{2}$-trace of the heat semi-group is obtained
by integrating the fibre trace ${\rm Tr}_{V}$ of the heat kernel
diagonal $U(x,x;t)$ over $M$, and reads as
\begin{equation}
{\rm Tr}_{L^{2}}\Bigr(e^{-t{\widetilde P}}\Bigr) =\int_{M}dx {\rm
Tr}_{V}U(x,x;t).
\end{equation}

In our paper, following \cite{gilk95b}, we are interested in a
slight generalization of the previous equation, where
$e^{-t{\widetilde P}}$ is `weighted' with a smooth scalar function
$f$ on $M$. More precisely, we are interested in the asymptotic
expansion as $t \rightarrow 0^{+}$ of the functional trace
\begin{equation}
{\rm Tr}_{L^{2}}\Bigr(f e^{-t{\widetilde P}}\Bigr) =\int_{M}dx  \;
f(x){\rm Tr}_{V} U(x,x;t).
\end{equation}
The results for the original problem are eventually recovered by
setting $f=1$, but it is crucial to keep $f$ arbitrary throughout
the whole set of calculations, as will be clear from the following
sections.

The asymptotic expansion of such a functional trace has the form
\begin{equation}
{\rm Tr}_{L^{2}}\Bigr(f e^{-t {\widetilde P}}\Bigr) \sim
\sum_{n=0}^{\infty}t^{(n-m)/2} a_{n}(f,{\widetilde P},{\mathcal
B}_\theta).
\end{equation}
Note that there is a change of convention in the indexing of
the Seeley coefficients with respect to the work in
\cite{espo02-66-085014}, i.e. our $a_{n}$ is denoted therein
by $a_{n/2}$.
The coefficients $a_{n}(f,{\widetilde P},{\mathcal B}_\theta)$
consist of two different parts, the interior part $a_{n}^{M}
(f,\widetilde P)$ and the boundary part $a_{n}^{\partial
M}(f,{\widetilde P},{\mathcal B}_\theta)$, i.e.
\begin{equation}
a_{n}(f,{\widetilde P},{\mathcal B}_\theta)
=a_{n}^{M}(f,{\widetilde P}) +a_{n}^{\partial M}(f,{\widetilde
P},{\mathcal B}_\theta).
\end{equation}
The interior parts $a_{n}^{M}(f,{\widetilde P})$ are obtained by
integrating some geometric invariants (see below) over $M$ and are
independent of the boundary conditions. By contrast, the boundary
parts $a_{n}^{\partial M}(f,{\widetilde P},{\mathcal B}_\theta)$
are obtained by integrating some geometric invariants over the
boundary $\partial M$ and these parts depend in a crucial way on
the boundary conditions. They will be the main concern of our
research from now on. The interior invariants are built
universally and polynomially from the metric tensor, its inverse,
the Riemann curvature of $M$, the bundle curvature (if a vector
bundle over $M$ is given) and the endomorphism (or `potential'
term) in the squared operator $P^{2}$. By virtue of Weyl's work on
the invariants of the orthogonal group, these polynomials can be
formed by using only tensor products and contraction of tensor
arguments. Here, the structure group is $O(m)$, $m$ being the
dimension of $M$. However, when a boundary occurs, the boundary
structure group is $O(m-1)$, and the Weyl theorem is used again to
construct invariants.

The structure of the article is as follows. In Section 2 we write
down the general form of the leading coefficients $a_1$ and $a_2$.
The special case calculation of \cite{espo02-66-085014} and
different functorial techniques are used to determine part of the
numerical multipliers of the geometric invariants. Further special
cases are shown in Section 3 and the complete $a_1$ and $a_2$
coefficients are determined. We end the paper with concluding
remarks.

\section{Determination of the leading coefficients}
We first write down the general form of the leading two boundary
contributions to the heat kernel (hereafter, $L_{aa}$ is our
notation for the trace of the extrinsic-curvature tensor of the
boundary).
\begin{lemma}\label{lem2.1}
Let $f$ be scalar. There exist universal constants $c_i(\theta,m)$
such that (hereafter our notation for the invariant integration
measure on $\partial M$ is simply $dy$) \beq a_1^{\partial M}
(f,\widetilde P,{\mathcal B}_\theta ) = (4\pi)^{-(m-1)/ 2}
\int\limits_{\partial M}dy \Tr_V ( c_1 (\theta
,m) f )  ,\label{coe1}\\
a_2^{\partial M} (f,\widetilde P,{\mathcal B}_\theta ) =
(4\pi)^{-m/2} \int\limits_{\partial M}dy \Tr_V \left( c_2 (\theta
, m) L_{aa} f +
c_3 (\theta , m) f \psi \tilde \gamma \gamma_m \right.\nn\\
\left. + c_4 (\theta , m ) f \psi \gamma_m+ c_5 (\theta , m) f
\psi \tilde \gamma + c_6 (\theta , m) f \psi + c_7 (\theta , m )
f_{;m} \right).\label{coe2}
\eeq
\end{lemma}
\begin{proof} This is a direct consequence of the
Weyl theorem on the invariants of the orthogonal group
\cite{gilk95b}, as we said at the end of Section 1.\end{proof} We
next determine the universal multipliers $c_i (\theta,m)$,
$i=1,...,7$. We first exploit a known special case. As usual, the
hypergeometric function is denoted by $_2F_1 (a,b;c;z)$.
\begin{lemma}\label{lem2.2}
We have
\beq
\quad\quad c_1 (\theta ,m) = \frac 1 4 \left( \cosh
^{m-1} \theta -1\right), \label{uni1}\\ \quad\quad c_2 (\theta ,
m) = \frac 1 {2(m-1)} \left\{ \frac{ 2m-5} 3 + (2-m) \,\, _2F_1
\left( 1 , \frac{m-1} 2 ; \frac 3 2 ; \tanh ^2 \theta \right)
\right\} .\label{uni2}
\eeq
\end{lemma}
\begin{proof} In \cite{espo02-66-085014} the heat kernel
coefficients for the given setting have been evaluated on the
Euclidean ball for the case $\psi =0$ and $f=1$. The results
obtained were \beq a_1 &=& \frac{ \sqrt \pi d_s}{2^m \Gamma
(m/2)} \left( \cosh^{m-1} \theta -1\right) , \nn\\
a_2 &=& \frac{(2m-5)d_s}{3\cdot 2^m \Gamma (m/2)} + \frac{
d_s}{2^m \Gamma (m/2)} \times\nn\\
& &\left\{ _2F_1 \left( 1,\frac{m-1} 2 ; \frac 1 2 ; \tanh ^2
\theta \right) - (m-1) \,\, _2F_1 \left( 1,\frac{m+1} 2 ; \frac 3
2 ; \tanh ^2 \theta \right) \right\}.\nn \eeq The volume of the
sphere, which is the boundary of the ball, is
$$
\mbox{vol} (S^{m-1}) = \frac{2\pi^{m/2}}{\Gamma(m/2)}.
$$
Using this to rewrite the coefficient $a_1$ shows the
assertion (\ref{uni1}), which agrees with Eq. (40) in \cite{anto90}
for the first boundary correction to the partition function.

To show (\ref{uni2}) we first use the Gauss recursion formula, see
e.g.~\cite{grad65b}, equation 9.137.17,
\beq
\gamma \,\, _2 F_1
(\alpha , \beta ; \gamma; z) - (\gamma -\beta  ) \,\, _2 F_1
(\alpha, \beta; \gamma +1; z) - \beta \,\, _2F_1 (\alpha , \beta
+1 ; \gamma +1 ; z) =0 ,\nn\eeq to write $a_2$ for the $m$-ball as
\beq a_2 = \frac{ d_s}{2^m \Gamma (m/2)} \left\{ \frac{2m-5} 3 +
(2-m) \,\, _2 F _1 \left( 1 , \frac{ m-1} 2 ; \frac 3 2 ; \tanh ^2
\theta \right) \right\} .\nn
\eeq
Comparison with  the general form
(\ref{coe2}) then shows assertion (\ref{uni2}). Note that in the
given setting, i.e. with $\psi =0 $ and $f=1$, the $c_2 (\theta
, m)L_{aa} f$ term is the only term contributing.
\end{proof}
\begin{remark}\label{rem2.3} For $\theta =0$
the boundary conditions reduce to
standard boundary conditions of mixed type. For $\theta =0$ we
have
\beq c_1 (0,m) &=&0,\nn\\
c_2 (0,m) &=& \frac 1 {2(m-1)} \left\{ \frac{ 2m-5} 3 + (2-m)
\cdot 1\right\} = \frac 1 {2(m-1)} \frac{1-m} 3 = -\frac 1 6
.\nn
\eeq
To achieve comparison with the known results for mixed
boundary conditions note that the auxiliary Hermitian endomorphism
$\chi$ needed to define the splitting of the spinor bundle is
\cite{kirs01b}
\beq
\chi = -\tilde \gamma \gamma_m.\nn
\eeq
Let
$$
\Pi _\pm = \frac 1 2 \left( 1 \pm \chi \right)
$$
be the projection on the $\pm$ eigenspaces of $\chi$. Mixed boundary
conditions are then defined as
\beq
{\mathcal B} \varphi = \Pi _-
\varphi \left| _{\partial M} \oplus ( \nabla_m + S) \Pi _+ \varphi
\right| _{\partial M} =0.\nn
\eeq
The relevant $S$ for the given setting is
$$
S = -\frac 1 2 L_{aa} \Pi _+ .
$$
Using the fact that $\Tr_V (\chi) =0$, $\Tr_V (\Pi _\pm ) = d_s
/2$, the coefficients for the relevant mixed boundary conditions
follow to be \beq a_1^{\partial M} (1, \widetilde P, {\mathcal
B}_0) &=&
0 , \nn\\
a_2 ^{\partial M} (1,\widetilde P,{\mathcal B}_0) &=&
(4\pi)^{-m/2} \frac 1 6 \int\limits_{\partial M} \Tr_V ( 2L_{aa} +
12 S ) dy =(4\pi)^{-m/2} \frac 1 6 \int\limits_{\partial M} \Tr_V
( -L_{aa})dy, \nn \eeq in agreement with our findings for $c_1
(0,m)$ and $c_2 (0,m)$.\end{remark} We next exploit the fact that
the connection $\nabla$ is not canonically defined. To simplify
the notation slightly we assume a localizing function $f=1$.
\begin{lemma}\label{lem2.4} We have
$$
c_4 (\theta , m) =0.
$$
\end{lemma}
\begin{proof} Let $\sigma_i$ be a skew-adjoint endomorphism of the
spinor bundle commuting with the Clifford structure $\gamma$,
$[\sigma_i , \gamma_j]=0$. Then $$\nabla_i (\epsilon ) = \nabla_i
+ \epsilon \sigma_i $$ defines a smooth one-parameter family of
compatible unitary connections. We define
$$
\psi (\epsilon ) :=
\psi - \epsilon \gamma_i \sigma_i
$$ to ensure that
$$
P(\epsilon )
= \gamma_i \nabla_i (\epsilon) + \psi (\epsilon ) = P
$$
is unaffected by the perturbation; the boundary condition also
remains unchanged. Therefore, the heat trace coefficient
(\ref{coe2}) remains unchanged. Using
$\tilde\gamma\gamma_i=-\gamma_i\tilde\gamma$ we evaluate the
variation $\delta = (d/d\epsilon)|_{\epsilon =0}$ of the single
terms for $\sigma_a=0$, $\sigma_m \neq 0$: \beq & & \delta \Tr_V
(\ctw L_{aa})
=0,\nn\\
& &\delta\Tr_V (\cth \psi \tilde \gamma \gamma_m )= -\Tr_V
(\cth\gamma_m\sigma_m \gt \gamma_m ) = \Tr_V (\cth\sigma_m \gt )
\nn\\ & &\quad \quad = -\Tr_V (\cth \sigma _m \gamma_m \gt
\gamma_m
) = - \Tr_V (\cth \sigma_m \gt ) =0 ,\nn\\
& &\delta\Tr_V (\cfo \psi \gamma_m ) = -\Tr_V (\cfo \gamma_m
\sigma_m
\gamma_m ) = \Tr_V (\cfo \sigma_m ), \nn\\
& &\delta \Tr_V (\cfi \psi \gt ) = -\Tr_V (\cfi \gamma_m \sigma_m
\gt ) = - \Tr_V (\cfi \sigma_m \gt )\nn\\ & &\quad \quad = \Tr_V
(\cfi \sigma_m \gamma_m \gt )
= \Tr_V (\cfi \gamma_m \sigma_m \gt ) =0 ,\nn\\
& &\delta \Tr_V (\csi \psi ) = - \Tr_V (\csi\gamma_m \sigma_m ) =
-\Tr_V
(\csi\gamma_m \sigma_m \gt\gt)\nn\\
& &\quad\quad = \Tr_V (\csi\gt \gamma_m \sigma_m \gt)
 = \Tr_V (\csi\gamma_m \sigma_m ) =0 .\nn
 \eeq
For the coefficient to remain unchanged we need
$\cfo =0$.

Considering $\sigma_a \neq 0$ and $\sigma_m =0$ does not produce
any new information.
\end{proof}
To find more information about the remaining unknown multipliers,
one might enlarge the setting and allow for an endomorphism-valued
$f$. However, apart from the fact that the number of invariants
goes up to $36$ and the calculation gets cumbersome, this does not
produce any relevant information for our problem and we do not
present further details.

Instead, we next exploit conformal rescaling techniques.
\begin{lemma}\label{lem2.7}
We have
$$
c_7 (\theta ,m) = -\frac{m-1}{m-2} \left(\ctw + \frac 1
6 \right).
$$
\end{lemma}
\begin{proof} Let $f$ be a smooth function with $f|_{\partial M}
=0$. Define $ds^2 (\epsilon ) := e^{2\epsilon f} ds^2$ and
$P(\epsilon ) := e^{-\epsilon f} P$. Let $\nabla$ be a compatible
unitary connection. We expand $P = \gamma^\nu
\nabla_{\partial_\nu} + \psi$ with respect to a local coordinate
system $x=(x_1,...,x_m)$ and use the metric to lower indices and
define $\gamma_\nu$. If we define
$$
\nabla (\epsilon )
_{\partial_\mu } := \nabla_{\partial _\mu} + \frac 1 2 \epsilon (
f_{;\nu} \gamma^\nu \gamma_\mu + f_{;\mu}),
$$
results of \cite{dowk99-242-107} show that $\nabla (\epsilon )$ is
a compatible unitary connection. Furthermore,
$$
\psi (\epsilon ) =
e^{-\epsilon f} \left( \psi - \frac 1 2 \epsilon (m-1) f_{;\nu}
\gamma^\nu \right).
$$
Note that the boundary condition remains unchanged under conformal
variation. The heat kernel coefficients satisfy the equation \beq
\left. \frac d {d\epsilon} \right|_{\epsilon =0} a_n \left( 1,
\widetilde P (\epsilon ) , {\mathcal B}_\theta \right) = (m-n) a_n
\left( f , \widetilde P , {\mathcal B}_\theta \right)
.\label{conf1} \eeq To study the numerical multiplier $c_7 (\theta
,m)$ we need the variations \beq \left. \frac d {d\epsilon}
\right|_{\epsilon =0} \tau (\epsilon ) &=& -2f\tau -
2(m-1) \Delta f ,\nn\\
\left. \frac d {d\epsilon} \right|_{\epsilon =0}L_{aa} (\epsilon )
&=& - f L_{aa} - (m-1) f_{;m} ,\nn \eeq where $\tau = R_{ijji}$ is
the scalar curvature. Applying equation (\ref{conf1}) shows the
assertion.
\end{proof}
\begin{remark} Note that, despite the appearance, the multiplier
$c_7 (\theta ,m)$ is well defined in dimension $m=2$. Using the
result for $\ctw$ given in equation (\ref{uni2}) we obtain
explicitly \beq c_7 (\theta ,m) = -\frac 1 2 \left\{ 1 - \,\,_2
F_1 \left( 1 , \frac{m-1} 2 ; \frac 3 2 ; \tanh ^2 \theta \right)
\right\}.\label{coec7} \eeq For $\theta =0$ this agrees with the
previous computation for mixed boundary conditions.
\end{remark}

\section{Relating the zeta and eta invariants}

In order to determine the numerical multipliers $\cth$, $\cfi$ and
$\csi$ we relate the zeta invariant to the eta invariant. We will
then evaluate the eta invariant on the $m$-dimensional cylinder
and ball for the case of an endomorphism-valued $f$. On the ball
we will restrict to the choices $f=1$ and $f=\gt$, respectively,
which will allow us to find $\cfi$ and $\csi$. Instead, on the
cylinder we can deal with general $f$. Performing the two special
case calculations is strictly speaking not necessary, but provides
helpful crosschecks of the answers obtained.

To distinguish the coefficients in the heat trace,
$\mbox{Tr}_{L^2} (f e^{-tP^2})$, and in the trace related to the
eta invariant, $\mbox{Tr}_{L^2} (fPe^{-tP^2})$, in this section we
use the notation \beq \mbox{Tr}_{L^2} \left( f e^{-tP^2}\right)
\sim \sum_n t^{(n-m)/2}
a_n^\zeta \left( f , P^2 , {\mathcal B}_\theta \right), \nn\\
\mbox{Tr}_{L^2} \left( f P e^{-tP^2}\right) \sim \sum_n
t^{(n-m-1)/2} a_n ^\eta \left( f , P , {\mathcal B} _\theta
\right) .\nn\eeq

The result we are going to need is the following:
\begin{lemma}\label{lem3.1}
Let $f \in C^\infty (\mbox{End} (V))$ and let $P_\epsilon := P +
\epsilon f$. We then have
\beq
\partial _\epsilon a_n ^\zeta (1, P^2_\epsilon , {\mathcal B}_\theta
) = - 2 a_{n-1} ^\eta (f , P_\epsilon , {\mathcal B}_\theta
).\nn
\eeq
\end{lemma}
\begin{proof}
The proof is insensitive to the boundary conditions imposed and
parallels the proof in \cite{bran92-108-47}.
\end{proof}
\begin{remark}\label{rem3.2} The very useful property of
this result is that the $a_n^\zeta$ coefficient for the zeta
invariant is related to the coefficient $a^\eta_{n-1}$ for  the
eta invariant, which will have a significantly simpler structure.
\end{remark} In order to apply Lemma \ref{lem3.1} to the
coefficient $a_2^\zeta$ we need the general form of the $a_1^\eta$
coefficient.
\begin{lemma}\label{lem3.3} Let $f\in C^\infty (\mbox{End} (V))$.
There exist universal constants $d_i (\theta , m)$ such that
\beq
& & a_1 ^{\eta , \partial M} (f,P,{\mathcal B}_\theta ) = \nn\\
& &\quad (4\pi)^{-m/2} \int\limits_{\partial M} dy \; \Tr_V
\left\{ \don f + \dtw f\gt + \dth f\gamma_m + \dfo f \gt \gamma_m
\right\} .\nn \eeq
\end{lemma}
\begin{proof} This follows immediately from the theory of
invariants taking into account that $f$ is in general a
matrix-valued endomorphism. \end{proof}
\begin{remark}\label{rem3.4} Lemma \ref{lem3.1} relates the
universal constant $d_j (\theta ,m)$, $j=1,...,4,$ with $c_i
(\theta ,m)$, $i=3,...,6$. In particular we have
\beq
& &\cth = -2
\dfo , \quad \cfo = - 2\dth , \nn\\
& & \cfi = -2 \dtw , \quad \csi = -2 \don .\nn \eeq From Lemma
\ref{lem2.4} we conclude $\dth = 0$. We evaluate $\don$ and $\dtw$
for the example of the ball and thus find $\cfi$ and $\csi$. We
also evaluate $\don$, $\dtw$ and $\dfo$ for the example of the
cylinder. This provides checks of the answers for $\cfi$ and
$\csi$ and in addition determines $\dfo$ and thus $\cth$.
\end{remark} For the case $f=1$ we proceed as described in
\cite{kirch}. The case $f=\gt$ is based upon this calculation and
therefore we need to present some details for the case $f=1$. We
first summarize properties of the spectral resolution for the
Dirac operator on the ball. Let $P =\gamma_i \nabla_i$ be the
Dirac operator on the ball and let us denote by $\varphi_{\pm}$
its eigenfunctions obeying the eigenvalue equation $P \varphi_\pm
= \pm \mu \varphi _\pm$. On writing the eigenvalue equation in
this form we have $\mu >0$. Later on we will write the eigenvalues
of $P$ as $\lambda = \pm \mu$, such that $|\lambda|=\mu$. Modulo a
suitable radial normalizing constant $C$, we may express
\cite{dowk96-13-2911}
\begin{eqnarray}
\varphi_{\pm}^{(+)}&=&{\frac{C}{r^{(m-2)/2}}} \left(
     \begin{array}{c}
        iJ_{n+m/2}(\mu r)
       \,Z^{(n)}_+(\Omega ) \\
     \pm J_{n+m/2-1}(\mu r)\,Z^{(n)}_+(\Omega )
        \end{array}  \right)  ,\text{ and} \label{sol1}\\
\varphi_{\pm}^{(-)}&=&{\frac{C}{r^{(m-2)/2}}}\left(
     \begin{array}{c}
     \pm J_{n+m/2-1}(\mu r)\,Z^{(n)}_
-(\Omega )  \\
   i J_{n+m/2}(\mu r)\,Z^{(n)}_-(\Omega ) \end{array}
\right). \label{sol2}
\end{eqnarray}
Here, $J_\nu (z)$ are the Bessel functions and $Z_\pm ^{(n)}
(\Omega )$ are the eigenspinors of the Dirac operator $\breve P$
on the sphere \cite{camp96-20-1},
\begin{eqnarray}
\breve P {\mathcal Z} _\pm ^{(n)} (\Omega )=  \pm \left(
n+\frac{m-1} 2 \right)
             {\mathcal Z} _\pm ^{(n)} (\Omega )\text{ for
             }n=0,1,...\nn
         \eeq
The degeneracy $d_n (m)$ for each eigenvalue is
\beq
d_n(m):=\dim
{\mathcal Z} _\pm ^{(n)} (\Omega )= \frac 1 2 d_s \left(
    \bea {c}
       m+n-2 \\
        n
      \eea \right) .\nn
\end{eqnarray}
We next apply the boundary operator which reads explicitly,
from Eq. (1.b),
\beq
\frac 1 2 \left(\begin{array}{cc} 1 & -ie^\theta \\
ie^{-\theta} & 1 \end{array}\right),\nn \eeq to the solutions
(\ref{sol1}) and (\ref{sol2}). This produces the following
eigenvalue conditions: \beq J_{n+\frac m 2 }(\mu ) \mp e^\theta
J_{n+\frac m 2 -1} (\mu ) &=& 0  \quad \mbox{for} \quad
\varphi_\pm
^{(+)} , \label{eigen1}\\
J_{n+\frac m 2 }(\mu ) \pm e^{-\theta} J_{n+\frac m 2 -1} (\mu )
&=& 0 \quad \mbox{for} \quad \varphi_\pm ^{(-)} . \label{eigen2}
\eeq These equations allow us to rewrite the eta function \beq
\eta (s;1, P,{\mathcal B}_\theta ) = \sum_\lambda \mbox{sgn}
(\lambda) |\lambda|^{-s} \nn \eeq in terms of a contour integral
and to apply the techniques described in detail in
\cite{bene02-35-9343,bord96-37-895,bord96-182-371,kirs01b}. The
coefficients in the asymptotic expansion for the eta invariant are
then determined by evaluating residues of $\eta$ according to
\cite{gilk95b} \beq \mbox{Res  }\eta (m-n;1, P, {\mathcal
B}_\theta ) = \frac{2 a_n^\eta (1, P, {\mathcal B}_\theta ) }
{\Gamma \left( \frac{ m-n+1} 2 \right) } .\label{etacoef} \eeq For
notational convenience we introduce $p=n+m/2 -1$. Starting point
of the analysis is \cite{kirch} \beq \eta (s; 1,P, {\mathcal
B}_\theta ) = \relsum \frac 1 {2\pi i} \int\limits_\Gamma dk
k^{-s} \frac d {dk} \ln \frac{1 + e^\theta \frac{ J_{p+1} (k)}{J_p
(k)} }{ 1 -e^\theta \frac{ J_{p+1} (k)}{J_p (k)}}- (\theta \to
-\theta), \nn \eeq where $\Gamma$ is a suitable counterclockwise
contour enclosing all solutions of the equations (\ref{eigen1})
and (\ref{eigen2}). After deforming the contour to the imaginary
axis this gives \beq \eta (s; P, {\mathcal B}_\theta )&=& \frac 1
{\pi i} \cos \left( \frac{ \pi s} 2 \right) \relsum
\int\limits_0^\infty dz z^{-s} \frac d {dz} \ln \frac{1-ie^\theta
\frac{I_{p+1} (z)}{I_p (z)}}
{1+ie^\theta \frac{I_{p+1} (z)}{I_p (z)}}- (\theta \to -\theta ) \nn\\
&&\hspace{-1.0cm}= \frac 1 {\pi i} \cos \left( \frac{ \pi s} 2
\right) \relsum \int\limits_0^\infty dz z^{-s} \frac d {dz} \ln
\frac{1+\frac p z i e^\theta -ie^\theta \frac{I_{p}' (z)}{I_p
(z)}} {1-\frac p z ie^\theta +ie^\theta \frac{I_{p}' (z)}{I_p
(z)}}- (\theta \to -\theta ),\nn
\eeq
where in the last step we have used
the recursion for the modified Bessel function \cite{grad65b}
\beq
I_{\nu +1} (z) = I_\nu ' (z) - \frac \nu z I_\nu (z)
.\label{recbes}
\eeq
In order to recover the coefficient $a_1^\eta$
we only need to consider the leading term in the uniform
$p\to\infty$ asymptotic expansion of the Bessel function
\cite{abra70b},
\beq
\frac{I_p ' (kp)}{I_p (kp)} \sim \frac{
(1+k^2)^{1/2}} k \left( 1 + {\mathcal O} \left( \frac 1 p \right)
\right) . \nn
\eeq
Hence we only need to find the residue of
\beq
A_0(s;1)=
\frac 1 {\pi i } \cos \left( \frac{ \pi s } 2 \right)
\relsum p^{-s} \int\limits_0^\infty dk k^{-s} \frac d {dk} \ln
\frac{ 1+\frac i k e^\theta - \frac{ ie^\theta \sqrt{ 1+k^2}} k }{
1-\frac i k e^\theta + \frac{ ie^\theta \sqrt{ 1+k^2}} k } -
(\theta \to -\theta )\nn
\eeq
at $s=m-1$.

We first observe that the summation over $n$ produces a multiple
of the Barnes zeta function \cite{barn03-19-374}, which is defined by
\beq
\zeta_{\mathcal B} (s,a) := \sum_{n=0}^\infty \left(
\begin{array}{c}
m+n-2 \\ n \end{array} \right) (n+a)^{-s} .\nn
\eeq
In detail we have
\beq
\relsum p^{-s} = \frac 1 2 d_s \zeta_{\mathcal B} \left(
s , \frac m 2 -1 \right).\nn
\eeq
In order to perform the
$k$-integral we first combine $\theta$ and $-\theta$ and evaluate
the logarithmic derivative to give
\beq
& &\frac d {dk} \ln \frac{
1+\frac{ie^\theta} k \left( 1-\sqrt{1+k^2}\right)} {
1-\frac{ie^\theta} k \left( 1-\sqrt{1+k^2}\right)}- (\theta \to
-\theta ) = -\frac{4i\sinh \theta}{\sqrt{1+k^2} \left( 2 + k^2 +
k^2 \cosh (2\theta)\right)} .\nn
\eeq
The relevant $k$-integral
therefore reads as
\beq
& &\int\limits_0^\infty dk k^{-s}
(1+k^2)^{-3/2} \frac 1 {1+\frac{k^2}{2(1+k^2)} \left( \cosh
(2\theta ) -1
\right)} \nn\\
&&\quad \quad =\frac 1 2 \Gamma \left( 1+\frac s 2 \right)
\sum_{l=0}^\infty (-1)^l \left( \frac{ \cosh (2\theta ) -1} 2
\right) ^l \frac{ \Gamma \left( \frac {1-s} 2 +l \right) } {\Gamma
\left( \frac 3 2 +l \right) } \nn\\
& &\quad \quad =\frac 1 {\sqrt{\pi}} \Gamma \left( 1+\frac s 2
\right) \Gamma \left( \frac{ 1-s} 2 \right) \,\, _2F_1 \left( 1,
\frac{1-s} 2 ; \frac 3 2 ; \frac 1 2 (1-\cosh (2\theta))\right)
\nn\\
& &\quad \quad =\frac 1 {\sqrt{\pi}} \Gamma \left( 1+\frac s 2
\right) \Gamma \left( \frac{ 1-s} 2 \right) \,\, _2F_1 \left( 1,
\frac{1-s} 2 ; \frac 3 2 ; -\sinh ^2 \theta \right) .\nn
\eeq
From here, with \cite{grad65b}
$$
\Gamma \left( \frac{ 1-s} 2 \right) =
\frac \pi {\cos \left( \frac{\pi s} 2 \right) \Gamma \left( \frac{
1+s} 2 \right)},
$$
we easily compute
\beq
A_ 0 (s;1) &=& - \frac 1
{\sqrt \pi} d_s \frac{ \Gamma \left( 1+\frac s 2 \right) } {
\Gamma \left( \frac {1+s} 2 \right)} \sinh \theta \,\, _2 F _1
\left(1 , \frac{ 1-s} 2 ; \frac 3 2 ; - \sinh ^2 \theta \right)
\zeta _{\mathcal B} \left( s , \frac m 2 -1 \right) . \nn
\eeq
On using \cite{grad65b}
\beq
\frac{\Gamma \left( \frac{m-1} 2
\right)} {\Gamma (m-1)} = \frac{ \sqrt \pi}{2^{m-2} \Gamma \left(
\frac m 2\right)} ,\nn
\eeq
the coefficient $a_1^\eta$ can be
cast into the form
\beq
a_1^\eta (1, P,{\mathcal B}_\theta ) &=&
\frac 1 2 \Gamma \left( \frac m 2\right) \mbox{Res } \eta (m-1;
1,P,{\mathcal B}_\theta ) = \frac 1 2 \Gamma
\left( \frac m 2 \right)  \mbox{Res } A_0 (m-1;1) \nn\\
&=& - \sinh \theta d_s \frac{ m-1} {2^m \Gamma \left( \frac m
2\right) } \,\, _2 F _1 \left( 1, 1-\frac m 2 ; \frac 3 2 ; -\sinh
^2 \theta \right) .\nn \eeq Comparing this with the answer on the
ball expected from Lemma \ref{lem3.3}, \beq a_1^\eta (1, P ,
{\mathcal B}_\theta ) = (4\pi ) ^{-m/2} \mbox{vol} (S^{m-1}) d_s
d_1 (\theta , m) = \frac 2 {2^m \Gamma \left( \frac m 2 \right)}
d_s d_1 (\theta ,m) , \nn \eeq we read of \beq d_1 (\theta ,m) =
-\frac{m-1} 2 \sinh \theta \,\, _2 F_1\left( 1, 1-\frac m 2 ;\frac
3 2 ; -\sinh ^2 \theta \right) .\label{d1ball} \eeq From Remark
\ref{rem3.4} we then get \beq \csi = (m-1) \sinh \theta \,\, _2 F
_1 \left(1, 1-\frac m 2 ; \frac 3 2 ; -\sinh ^2 \theta \right)
.\nn \eeq

To find the universal constant $\cfi$ we perform the calculation
on the ball with $f=\gt$. This choice complicates the analysis
significantly because the normalization constant $C$ and further
integrals over products of Bessel functions come into the play.
First we note that if $\eta (s; x,y)$ denotes the local eta
function, then \beq \eta (s; x,y) = \sum_{\mu}  \mu ^{-s} \left\{
\varphi_+ ^{(\pm )} (x) ^* \,\, \varphi_+ ^{(\pm )} (y) -
\varphi_- ^{(\pm )} (x) ^* \,\, \varphi_- ^{(\pm )} (y)\right\}
.\nn \eeq We want to analyze \beq \Tr \left( \gt \eta (s;
x,x)\right) = \sum_{\mu} \mu ^{-s} \left\{< \varphi_+ ^{(\pm )} |
\gt \varphi_+ ^{(\pm )}
> -< \varphi_- ^{(\pm )} | \gt \varphi_- ^{(\pm )}
> \right\} ,\label{treta}
\eeq
with $<\varphi_1|\varphi_2>$
denoting the Hilbert space product
\beq
< \varphi _1 | \varphi _2>
\equiv \int\limits_M dx \; \varphi _1 ^* (x) \varphi_2 (x) . \nn
\eeq
Since
\beq
\gt =\left(
\begin{array}{cc} 1 & 0\\ 0 & -1 \end{array} \right) \nn
\eeq
changes the sign of the lower chirality, the normalization
constant $C$ does not cancel in the Hilbert space products
appearing in (\ref{treta}), but instead values of various
integrals occur explicitly. We first observe that
$$
\frac 1 {C^2} = \int\limits_0^1 dr \,\, r (J_{p+1} ^2 (\mu r) +
J_p ^2 (\mu r) ) .
$$
We use \cite{grad65b}
$$
\int\limits_0^1 dr \,\, r J_\nu ^2 (\mu r) = \frac 1 2 \left\{
J_\nu ^2 (\mu ) - J_{\nu -1} (\mu ) J_{\nu +1} (\mu ) \right\} $$
to find $$\frac 1 {C^2} = \frac 1 2 \left\{ J_p ^2 (\mu ) +
J_{p+1} ^2 (\mu ) - J_{p-1}  (\mu )J_{p+1} (\mu ) - J_p (\mu )
J_{p+2} (\mu ) \right\} .
$$
We use the implicit eigenvalue
equations (\ref{eigen1}) and (\ref{eigen2}) together with
recursion relations for the Bessel functions \cite{grad65b}
$$
J_{p+2} (\mu ) = \frac{ 2(p+1)} \mu J_{p+1} (\mu ) - J_p (\mu ), \quad
J_{p-1} (\mu ) = \frac{ 2p} \mu J_p (\mu ) - J_{ p+1} (\mu ) ,
$$
to simplify the normalization constants $C_{\pm}^{(\pm )}$ for the
different spinors $\varphi_\pm ^{(\pm )}$. We obtain \beq C_+
^{(\pm )} &=& \frac {\sqrt \mu}{J_p (\mu )} \frac 1 {\left(\mu+\mu
e^{\pm 2\theta} \mp (2p+1)
e^{\pm \theta }\right)^{1/2}} , \nn\\
C_- ^{(\pm )} &=& \frac {\sqrt \mu}{J_p (\mu )} \frac 1 {\left(\mu
+\mu e^{\pm 2\theta} \pm (2p+1) e^{\pm \theta }\right)^{1/2}} .
\nn \eeq Proceeding in the same way for the quantities
$<\varphi_\pm ^{(\pm )} | \gt \varphi_\pm ^{(\pm)}>$, we find \beq
<\varphi_+ ^{(\pm )} | \gt \varphi_+^{(\pm )} > &=& - \frac{
e^{\pm \theta}}{\mu+\mu e^{\pm 2\theta} \mp (2p+1) e^{\pm \theta}}
= - \frac 1 {2\cosh \theta} \frac 1
{\mu\mp \frac{ p+1/2} {\cosh \theta} } , \nn\\
<\varphi_- ^{(\pm )} | \gt \varphi_-^{(\pm )} > &=&  \frac 1
{2\cosh \theta} \frac 1 {\mu\pm \frac{ p+1/2} {\cosh \theta} } .
\nn \eeq Using these results in (\ref{treta}), we obtain the
following contour integral representation: \beq \eta ( s; \gt , P
, {\mathcal B}_\theta ) &=& - \frac 1 {4\pi i\cosh \theta }
\relsum \int\limits_\Gamma dk \,\, k^{-s} \frac{ \frac d {dk} \ln
\left[ J_{p+1} (k ) - e^\theta J_p (k ) \right]} {k-\frac
{p+1/2}{\cosh
\theta} }  \nn\\
& &- \frac 1 {4\pi i\cosh \theta } \relsum \int\limits_\Gamma dk
\,\, k^{-s} \frac{ \frac d {dk} \ln \left[ J_{p+1} (k ) + e^\theta
J_p (k ) \right]} {k+\frac {p+1/2}{\cosh \theta} }
\label{etagt}\\
& &+ (\theta \to -\theta ).\nn
\eeq
Note that the counterclockwise
contour must only include the zeroes of the equations
(\ref{eigen1}) and (\ref{eigen2}) such that the appropriate
summation over eigenvalues results. The poles at $k=(p+1/2) /
\cosh\theta$ should lie {\it outside} the contour because they have
been introduced by the normalization integral and need {\it not}
be summed over. The situation is similar to the analysis for
radial smearing functions, see \cite{dowk01-42-434} for more
details. This observation is important because when shifting the
contour towards the imaginary axis {\it additional} contributions
result. Using the index $p$ for all Bessel functions an
intermediate result reads as
\beq
\eta (s; \gt , P , {\mathcal
B}_\theta ) &=& \nn\\
& &\hspace{-2.5cm}\frac 1 {2\pi i \cosh \theta } \cos \left(
\frac{ \pi s} 2 \right) \relsum \int\limits_0^\infty dz \,\,
z^{-s} \frac{ \frac d {dz} \ln \left[ I_p ' (z) - \frac p z I_p
(z) - ie^\theta I_p (z) \right]} {iz + \frac{ p+1/2} {\cosh
\theta}}
\nn\\
& &\hspace{-2.5cm}+\frac 1 {2\pi i \cosh \theta } \cos \left(
\frac{ \pi s} 2 \right) \relsum \int\limits_0^\infty dz \,\,
z^{-s} \frac{ \frac d {dz} \ln \left[ I_p ' (z) - \frac p z I_p
(z) + ie^\theta I_p (z)
\right]} {iz - \frac{ p+1/2} {\cosh \theta} }\nn\\
& &\hspace{-2.5cm}+ \frac 1 {2 \cosh \theta }\relsum
\left(\frac{(p+1/2)}{\cosh \theta}\right)^{-s} \frac d {dk} \ln
\left.\left[J_p ' (k) + \left( e^\theta -\frac p k \right)J_p
(k)\right]\right|_{k=\frac{ p+1/2}{\cosh \theta}}\nn\\ & &
\hspace{-1.5cm}+ (\theta \to -\theta ).\label{etagt1}
\eeq
The last
contribution resulting from the shifting of the contour can be
given in closed form by using the differential equation for the
Bessel function \cite{grad65b},
$$
\left[{d^{2}\over dz^{2}}+{1\over z}{d\over dz}
+\left(1-{\nu^{2}\over z^{2}}\right)\right]J_{\nu}(z)=0.
$$
We calculate
\beq
& &\frac d {dk} \ln \left.\left( J_p ' (k) + \left(e^\theta -
\frac p k \right) J_p (k) \right) \right|_{k=\frac{ p+1/2} {\cosh
\theta}} \nn\\
& &\hspace{1.0cm} =\left.\frac{ J_p '' (k) + \frac p {k^2} J_p (k)
+ \left( e^{\theta} - \frac p k \right) J_p ' (k) }{J_p ' (k) +
\left(e^\theta - \frac p k \right) J_p (k)
}\right|_{k=\frac{p+1/2}
{\cosh \theta}} \nn\\
& &\hspace{1.0cm}=\left.\frac{ J_p ' (k) \left( e^\theta - \frac{
p+1} k \right) + J_p (k) \left( \frac{ p (p+1)}{k^2} -1 \right)} {
J_p ' (k) + \left( e^\theta - \frac p k \right) J_p (k)
}\right|_{k=\frac{ p+1/2} {\cosh \theta}} = \sinh \theta -
\frac{\cosh \theta} {2p+1} .\nn
\eeq
Adding the contributions from
$\theta $ and $-\theta$ the $\sinh \theta$ terms cancel and the
summation over $n$ leads to $\zeta_{\mathcal B} (s+1,(m-1)/2)$,
which has no pole at $s=m-1$. Therefore, for the present purpose
this term is irrelevant.

In the remaining integrals in (\ref{etagt1}) we need, as before,
only the leading term in the Debye asymptotic expansion of Bessel
functions. Explicitly, with $x=1/\cosh \theta$, we obtain to
leading order
\beq
A_0 (s;\gt ) &=& - \frac{ \cos \left( \frac{
\pi s} 2 \right)} {\pi \cosh \theta} \relsum p^{-s}
\int\limits_0^\infty dk \,\, k^{-s-1} \sqrt{ 1+k^2} \left[ \frac 1
{k-ix} + \frac 1
{k+ix}\right] \nn\\
&=& -\frac{ 2\cos \left( \frac{ \pi s} 2 \right)} {\pi \cosh
\theta }\relsum p^{-s} \int\limits_0^\infty dk \,\, \frac{k^{-s}
\sqrt{ 1+k^2}} {k^2 + x^2} .\nn \eeq The $k$-integral is
\cite{grad65b} \beq \int\limits_0^\infty dk \,\, \frac{ k^{-s}
\sqrt{1+k^2}}{k^2 +x^2} = \frac 1 {2x^2} \frac { \Gamma \left(
\frac{ 1-s} 2 \right) \Gamma \left( \frac s 2 \right)}{\Gamma
\left( \frac 1 2 \right)} \,\, _2F_1 \left( 1, \frac{ 1-s} 2 ;
\frac 1 2 ; 1-\frac 1 {x^2} \right) , \nn \eeq and hence \beq A_0
(s;\gt)= - \frac 1 2 d_s\frac{ \cosh \theta } {\sqrt \pi } \frac
{\Gamma \left( \frac s 2 \right)} {\Gamma \left( \frac {1+s} 2
\right)} \,\, _2 F_1 \left( 1 , \frac{ 1-s} 2 ; \frac 1 2 ; -\sinh
^2 \theta \right) \zeta _{\mathcal B} \left( s , \frac m 2 -1
\right).\nn \eeq The residue is easily evaluated and via
(\ref{etacoef}) we compare it with the form given in Lemma
\ref{lem3.3} to read off \beq d_2 (\theta ,m) = -\frac 1 2 \cosh
\theta \,\, _2 F _1 \left( 1 , 1-\frac m 2 ; \frac 1 2 ; -\sinh ^2
\theta \right),\label{d2ball} \eeq which implies \beq \cfi =\cosh
\theta \,\, _2 F _1 \left( 1 , 1-\frac m 2 ; \frac 1 2 ; -\sinh ^2
\theta \right). \nn \eeq
\begin{remark}\label{rem3.5} In the calculation just described it
is the argument $-\sinh ^2 \theta$ that occurs naturally in the
hypergeometric functions. Instead, the constant $\ctw$ in
(\ref{uni2}) and $\cse$ in (\ref{coec7}) have been given using
$\tanh^2 \theta$. In order to provide answers in a unified way one
might use the transformation formula \cite{grad65b}
$$
_2F_1 (\alpha , \beta ; \gamma ; z ) = (1-z)^{-\alpha} \,\, _2
F_1 \left(\alpha , \gamma -\beta ; \gamma ; \frac z {z-1}
\right)
$$
to write
\beq
\ctw &=& \frac 1 {2(m-1)} \left\{ \frac{
2m-5} 3 + (2-m) \cosh ^2 \theta \,\, _2F_1 \left( 1 , 2-\frac m 2
; \frac 3 2 ; -\sinh ^2 \theta \right) \right\},\nn\\
\cse &=& -\frac 1 2 \left \{ 1-\cosh ^2 \theta \,\, _2 F_1 \left( 1
, 2-\frac m 2 ; \frac 3 2 ; -\sinh ^2 \theta \right)\right \}.
\nn
\eeq
\end{remark}
\begin{remark} \label{rem3.6} Note that, despite the complicated
appearance of the universal constants, for each specific dimension
$m$ a simple function of $m$ and $\theta$ results. In particular,
whenever the second argument of $_2F_1$ is $0$ or a negative
integer, the hypergeometric function reduces to
a finite polynomial in $\sinh ^2\theta$.
\end{remark}

In order to find the missing multiplier $\cth$ we present a
calculation on the cylinder. In order to summarize previous
results \cite{bene03-36-11533} we need to provide some notation.
Let $M=\reals_+ \times N$ be an even dimensional cylinder equipped
with the metric $ds^2 = dx_m^2 + ds_N^2$, where $x_m$ is the
coordinate in $\reals_+$ and plays the role of the normal
coordinate, and $ds_N^2$ is the metric of the closed boundary $N$.
The coordinates on $N$ are denoted by $y=(y_1,y_2,...,y_{m-1})$.
To write down the heat kernel on $M$ for $P^2= (\gamma_i
\nabla_i)^2$ with boundary condition ${\mathcal B}_\theta$, we
call $\phi_\omega (y)$ the eigenspinors of the operator $B=\gt
\gamma_m \gamma_a \nabla _a$, corresponding to the eigenvalue
$\omega$, normalized so that \beq \sum_{\omega}
\phi_{\omega}^{\star}(y) \phi_{\omega}(y^{'})=
\delta^{m-1}(y-y^{'}), \nn\eeq with $\delta^{m-1}$ the Dirac delta
function, and \beq \int\limits_{N}\,\,dy\,\,
\phi_{\omega}^{\star}(y) \phi_{\omega}(y)=1\,.\nn\eeq Finally we
need $x=(y,x_m)$, $\xi = x_m - x_m '$, $\eta = x_m + x_m'$,
$u_\omega (\eta ,t) = \frac{\eta}{\sqrt {4t}} - \sqrt t \omega
\tanh \theta$, and the complementary error function
\[\mbox{erfc}(x)=\frac{2}{\sqrt{\pi}}\int_{x}^{\infty}d\xi
e^{-\xi^2}\, .\] We then have \cite{bene03-36-11533} \beq
U(x,x';t) &=& \frac{1}{\sqrt{4\pi
t}}\sum_{\omega}\phi_{\omega}^{\star}(y^{'}) \phi_{\omega}(y)
e^{-\omega^{2}
t}\left\{\left(e^{\frac{-\xi^2}{4t}}
-e^{\frac{-\eta^2}{4t}}\right)\mathbf{1}\right.\nn\\
& & \left.+\frac{2\pip \pipl}{\cosh^2(\Th)}\left[1+\sqrt{(\pi
t)}\omega \tanh{\Th}e^{u_{\omega}^{2}(\eta,t)}
\mbox{erfc}(u_{\omega}(\eta,t))\right]e^{\frac{-\eta^2}{4t}}\right\}
. \label{hkcyl}\eeq (Note that although the formal appearance of
the heat kernel is identical to the one in \cite{bene03-36-11533},
equation (5.1), the meaning of $\Pi_+$ is slightly different. The
reason is that \cite{bene03-36-11533} considers the boundary
condition resulting from, in our present notation, $\Pi_+$ whereas
we consider the one resulting from $\Pi_-$. Formally the
transition is obtained by reversing the sign of the normal and by
using our present notation for $\Pi_\pm$.) As remarked in
\cite{bene03-36-11533} the first term is the heat kernel on the
manifold $\reals\times N$, which does not encode any information
about the boundary contribution. In the following, without
changing the notation, we will ignore this term and we will
determine the boundary contributions to the eta invariant from the
remaining terms.

Let $f\in C^\infty (\mbox{End}(V))$, then we want to consider
$\mbox{Tr}_{L^2} (f [ P_x U (x,x'; t)]_{x=x'} )$; note that the
derivatives need to be performed before the coincidence limit
$x=x'$ is taken. Given we have the {\it local} form of the heat
kernel we can in principle deal with an arbitrary $f$. For our
present purpose it is easiest to assume $f=f(y)$ only such that
the $x_m$-integration can be done without complication.

It is natural to introduce the heat kernel $U_B (y,y'; t)$ of the
operator $B^2$, \beq U_B (y,y'; t) = \sum_\omega \phi _\omega ^*
(y') \phi _\omega (y) e^{-\omega ^2 t } ; \nn\eeq furthermore, to
make the single steps easier to follow we use the splitting
\beq
U_1 (x,x'; t) &=& -\frac{1}{\sqrt{4\pi
t}}\sum_{\omega}\phi_{\omega}^{\star}(y^{'}) \phi_{\omega}(y)
e^{-\omega^2 t}e^{\frac{-\eta^2}{4t}}, \nn\\
U_2 (x,x' ;t) &=& \frac{1}{\sqrt{4\pi
t}}\sum_{\omega}\phi_{\omega}^{\star}(y^{'}) \phi_{\omega}(y)
e^{-\omega^2 t}\frac{2\pip
\pipl}{\cosh^2(\Th)}\nn\\
& &\left[1+\sqrt{(\pi t)}\omega \tanh{\Th}e^{u_{\omega}^{2}(\eta,t)}
\mbox{erfc}(u_{\omega}(\eta,t))\right]e^{\frac{-\eta^2}{4t}} .
\nn\eeq Acting with $P$ and performing the $x_m$-integration,
intermediate results are \beq \int\limits_0^\infty dx_m \,\, f [
P_x U_1 (y,y', x_m, x_m'; t) ]_{x_m=x_m'} &=& \nn\\
& &\hspace{-5.0cm}\frac 1 {\sqrt {4\pi t}} \frac 1 2 f\gamma_m U_B
(y,y'; t) - \frac 1 4 f \gamma_m \gt
B_y U_B (y,y'; t) , \label{U1ex}\\
\int\limits_0^\infty dx_m \,\, f [ P_x U_2 (y,y', x_m, x_m'; t)
]_{x_m=x_m'} &=& \nn\\
& &\hspace{-7.0cm}-\frac 1 {2\cosh ^2 \theta } f \gamma_m U_B
(y,y';t) \Pi _+ \Pi _+^* \left[ \frac 1 {\sqrt{\pi t}} + \omega
\tanh \theta e^{t\omega^2 \tanh ^2\theta} \mbox{erfc}(-\sqrt t
\omega \tanh \theta )\right]\nn\\
& & \hspace{-5.0cm}+\frac 1 {2\cosh ^2 \theta } f \gamma_m \gt U_B
(y,y';t) \Pi _+ \Pi _+^* e^{t\omega^2 \tanh ^2\theta}
\mbox{erfc}(-\sqrt t \omega \tanh \theta ) . \label{U2ex}\eeq
Here, we have used the relation
\[ -\frac12 \frac{\partial}{\partial x_m}
\left[e^{-x_m^2/t+u_{\omega}^2(2x_m,t)} \mbox{erfc}
(u_{\omega}(2x_m,t))\right] =\]\beq
e^{-x_m^2/t}\left[\frac{1}{\sqrt{\pi t}}+\omega \tanh\theta\,
e^{u_{\omega}^2(2x_m,t)}
\mbox{erfc}(u_{\omega}(2x_m,t))\right]\,.\nn\eeq Whereas the
asymptotic $t\to 0$ behaviour in (\ref{U1ex}) could be easily
found from the corresponding (known) behaviour of the trace of
$U_B$, the same is not as simple for the result in (\ref{U2ex}).
We have found it most convenient to perform the $L^2 (N)$-trace
and to relate the above equations to the zeta and eta function via
\beq \zeta (s;f,P^2,{\mathcal B}_\theta ) &=& \mbox{Tr}_{L^2} ( f
(P^2 )^{-s}) = \frac 1 {\Gamma (s)} \int\limits_0^\infty dt \,\,
t^{s-1} \mbox{Tr}_{L^2} \left( f
e^{-t P^2}\right) , \nn \\
\eta (s;f,P,{\mathcal B}_\theta ) &=& \mbox{Tr}_{L^2} ( f P(P^2
)^{-s}) = \frac 1 {\Gamma \left( \frac{s+1} 2 \right)}
\int\limits_0^\infty dt \,\, t^{\frac{s-1} 2 } \mbox{Tr} _{L^2}
\left( f P e^{-tP^2} \right) ,\nn \eeq and to evaluate the
asymptotic $t\to 0$ expansion from (\ref{etacoef}) and \beq
\mbox{Res } \zeta (z; f,B^2 ) = \frac{ a_{\frac{m-1} 2 -z }^\zeta
(f , B^2)} {\Gamma (z)} .\label{zetaheatrel}\eeq For (\ref{U1ex})
the associated relation is readily found, \beq \eta _1
(s;f,P,{\mathcal B}_\theta ) = \frac 1 {\sqrt {4\pi}} \frac 1 2
\frac{\Gamma \left( \frac s 2 \right) } {\Gamma \left( \frac{ s+1}
2 \right)} \zeta \left( \frac s 2 ; f\gamma_m , B^2 \right)- \frac
1 4 \eta ( s ; f\gamma_m \gt , B ) .\nn\eeq In order to proceed
with (\ref{U2ex}) we note first that \beq t\omega^{2} \tanh ^2
\theta - t\omega^{2} &=& - \frac{
t\omega^2} {\cosh ^2 \theta } , \nn\\
\mbox{erfc} ( \sqrt t \omega \tanh \theta ) &=& 1 + \mbox{erf}
(\sqrt t \omega \tanh \theta ) , \nn\eeq
with the error function
\beq \mbox{erf} (x) = \frac 2 {\sqrt \pi} \int\limits_0 ^x dt \,\,
e^{-t^2} .\nn\eeq
The resulting $t$-integral then is
\beq
\int\limits_0^\infty dt \,\, t^{\frac{s-1} 2} e^{-\frac{t\omega^2}
{\cosh ^2 \theta} } \left( 1 + \mbox{erf} (\sqrt t \omega \tanh
\theta \right) = \nn\\\frac{ \cosh ^{s+1} \theta }{| \omega|
^{s+1}} \left[\Gamma \left( \frac{ s+1} 2 \right) + \frac 2 {\sqrt
\pi} \Gamma \left( 1 + \frac s 2 \right) \sinh \theta
\mbox{sgn}(\omega ) \,\, _2 F _1 \left( \frac 1 2 , 1 + \frac s 2
; \frac 3 2 ; -\sinh ^2 \theta \right) \right] .\nn
\eeq
This produces the following contributions to the eta function:
\beq
\eta _2 (s; f , P, {\mathcal B}_\theta ) = - \frac{ \Gamma \left(
\frac s 2 \right)} {2\sqrt \pi \Gamma \left( \frac{ s+1} 2 \right)
\cosh ^2 \theta } \zeta \left( \frac s 2 ; \Pi_+ \Pi _+ ^ * f
\gamma_m ,
B^2 \right) \nn\\
- \frac 1 2 \sinh \theta \cosh ^{s-2} \theta \eta ( s ; \Pi _+ \Pi
_+ ^* f \gamma_m , B ) \nn\\
- \frac 1 {\sqrt \pi} \frac{ \Gamma \left( 1 + \frac s 2
\right)}{\Gamma \left( \frac { s+1} 2 \right)} \sinh ^2 \theta
\cosh ^{s-2} \theta \,\, _2 F _1 \left( \frac 1 2 , 1+\frac s 2 ;
\frac 3 2 ; -\sinh ^2 \theta \right) \zeta \left( \frac s 2 ; \Pi
_+ \Pi _+ ^* f \gamma _m , B^2 \right) \nn\\
+ \frac 1 2 \cosh ^{s-1} \theta \eta ( s ; \Pi _+ \Pi _+ ^* f
\gamma _m \gt ; B ) \nn\\
+ \frac{ \cosh ^{s-1} \theta \sinh \theta \Gamma \left( 1 + \frac
s 2 \right)} { \sqrt \pi \Gamma \left( \frac {s+1} 2 \right)} _2
F_1 \left( \frac 1 2 ; 1 + \frac s 2 ; \frac 3 2 ; -\sinh ^2
\theta \right) \zeta \left( \frac s 2 ; \Pi _+ \Pi _+ ^* f
\gamma_m \gt , B^2 \right) .\nn\eeq From here, with the help of
(\ref{zetaheatrel}) and (\ref{etacoef}), it is easy to find the
residue of $\eta (s;f,P,{\mathcal B}_\theta )$ at $s=m-1$, needed
for the evaluation of $a_1^\eta ( f,P,{\mathcal B}_\theta )$. We
find \beq \mbox{Res } \eta (m-1; f , P , {\mathcal B}_\theta )
&=&\frac 1 {\sqrt \pi \Gamma \left( \frac m 2 \right)} \left\{
\frac 1 2 a_0 (f \gamma_m, B^2) - \frac 1 {\cosh ^2 \theta} a_0 (
\Pi _+ \Pi
_+^* f \gamma_m , B^2) \right.\nn\\
& &\hspace{-2.0cm}+(m-1) \sinh \theta \cosh ^{m-2} \theta \,\,
_2F_1 \left( \frac
1 2 , \frac{ m+1} 2 ; \frac 3 2 ; -\sinh ^2 \theta \right) \times\nn\\
& & \left.\left[ a_0 ( \Pi _+ \Pi _+ ^* f \gamma _m \gt , B^2 ) -
\tanh \theta a_0( \Pi _+ \Pi _+ ^* f \gamma _m  , B^2
)\right]\right\}.\nn\eeq The leading heat kernel coefficient $a_0
(G, B^2)$ is of course known for an arbitrary endomorphism $G$; it
is \beq a_0 (G,B^2) = (4\pi )^{-\frac{m-1} 2} \int\limits_N dy
\,\, \mbox{Tr}_V (G).\nn\eeq In order to obtain the invariant form
given in Lemma \ref{lem3.3}, we evaluate $\Pi_+\Pi_+^*$ in the
form
$$\Pi_+ \Pi_+^* = \frac 1 2 \cosh \theta (\cosh \theta + \gt \sinh
\theta - \gt \gamma_m ) .$$ Adding up all pieces this shows \beq
\mbox{Res }\eta (m-1;f,P,{\mathcal B}_\theta ) &=& \frac 1 {\sqrt
\pi \Gamma \left( \frac m 2 \right)} (4\pi )^{-m/2}
\int\limits_N dy \,\, \mbox{Tr}_V \left\{f \gamma_m \cdot 0 \right.\nn\\
& &\hspace{-4.0cm}+ f \gamma_m \gt \left[ \frac 1 2 \tanh\theta -
\frac 1 2 (m-1) \sinh \theta \cosh ^{m-2} \theta \,\, _2F_1 \left(
\frac 1 2 ,
\frac{m+1} 2 ; \frac 3 2 ; -\sinh ^2 \theta \right) \right]\nn\\
& &\hspace{-4.0cm}-f \frac 1 2 (m-1) \sinh \theta \cosh ^{m-1}
\theta \,\, _2 F _1\left( \frac 1 2 , \frac{m+1} 2 ; \frac 3 2 ;
-\sinh ^2 \theta
\right) \nn\\
& & \left.\hspace{-4.0cm} - f \gt \left[ \frac 1 {2\cosh \theta} +
\frac 1 2 (m-1) \sinh ^2 \theta \cosh ^{m-2} \theta \,\, _2 F
_1\left( \frac 1 2 , \frac{m+1} 2 ; \frac 3 2 ; -\sinh ^2 \theta
\right)\right]\right\}.\nn\eeq From this result we can read off
$d_i (\theta ,m)$, $i=1,...,4;$ we find \beq \don &=& - \frac
{m-1} 2 \sinh \theta \cosh^{m-1} \theta \,\, _2 F _1 \left( \frac
1
2 , \frac{m+1} 2 ; \frac 3 2 ; -\sinh ^2 \theta \right) , \nn\\
\dtw &=& - \frac 1 {2\cosh \theta} - \frac {m-1} 2 \sinh ^2 \theta
\cosh ^{m-2} \theta \,\, _2 F _1 \left( \frac 1 2 ,
\frac{m+1} 2 ; \frac 3 2 ; -\sinh ^2 \theta \right) , \nn\\
\dth &=& 0 , \nn\\
\dfo &=& - \frac 1 2 \tanh \theta + \frac {m-1} 2 \sinh \theta
\cosh ^{m-2} \theta \,\, _2 F _1 \left( \frac 1 2 , \frac{m+1} 2 ;
\frac 3 2 ; -\sinh ^2 \theta \right).  \nn\eeq The result for $d_1
(\theta ,m)$ can be seen to agree with the result on the ball, eq.
(\ref{d1ball}), by using the transformation formula
(\cite{grad65b}, eq. 9.131.1) \beq_2 F _1 ( \alpha , \beta ;
\gamma ; z) = (1-z )^{\gamma -\alpha -\beta }  \,\, _2 F_1 (
\gamma - \alpha , \gamma -\beta ; \gamma ; z)
.\label{transhyper}\eeq In order to show that the results for
$\dtw$ coming from the ball and cylinder agree, we need to show
that \beq\cosh ^2 \theta \,\, _2 F _1 \left( 1 , 1 -\frac m 2 ;
\frac 1 2 ; -\sinh ^2 \theta \right) \nn\\
= 1 + (m-1) \sinh ^2 \theta \cosh ^{m-1} \theta \,\, _2 F _1
\left( \frac 1 2 , \frac{m+1} 2 ; \frac 3 2 ; -\sinh ^2 \theta
\right) .\eeq To see this, we first apply the above transformation
formula, eq. (\ref{transhyper}), and then the Gauss recursion
formula (\cite{grad65b}, eq. 9.137.12)
$$ \gamma \,\, _2 F _1 ( \alpha , \beta ; \gamma ; z ) - \gamma
\,\, _2 F _1 (\alpha +1 , \beta ; \gamma ; z ) + \beta z \,\, _2 F
_1 (\alpha +1 , \beta +1 ; \gamma +1 ; z ) =0 $$ with $\alpha =
-1/2$, $\beta = (m-1)/2$, $\gamma = 1/2$, and $z=-\sinh^2\theta$.
Thus, all results obtained are consistent and we have determined
the full $a_1$ and $a_2$ coefficient for chiral bag boundary
conditions.

\section{Concluding remarks}

For the case of operators $P$ of Dirac type subject to local
boundary conditions of chiral bag type as in eq. (1.b), we have
studied the asymptotic expansion as $t \rightarrow 0^{+}$ of the
smeared $L^{2}$-trace of the associated heat semigroup, i.e.
$$
{\rm Tr}_{L^{2}}\Bigr(f e^{-tP^{2}}\Bigr) \sim \sum_{n=0}^{\infty}
t^{(n-m)/2}a_{n}(f,P^{2},{\mathcal B}_\theta).
$$
On using functorial methods, special case calculations and the
relation between $\eta$- and $\zeta$-invariants, we have succeeded
in evaluating the full boundary contribution to the $a_1$ and
$a_{2}$ coefficients, the functional form of which is given by
eqs. (\ref{coe1}) and (2.b). Our contributions are of technical
but non-trivial nature, because both functorial methods and the
theory of the $\eta$-invariant require a lot of work to obtain the
desired $a_{2}$ coefficient. It now appears possible that, by
exploiting the methods described in our paper, further heat-kernel
coefficients will be obtained, if they are needed in physical or
mathematical applications. In turn, a better understanding of the
spectral functions of modern mathematical physics \cite{kirs01b}
will also be gained.\\[.3cm]

{\bf\large Acknowledgements:} Research of GE was partially
supported by the INFN (Naples, Italy) and by PRIN 2002 {\it SINTESI}.
Research of PG was
partially supported by the MPI (Leipzig, Germany). KK acknowledges
support by the Baylor University Summer Sabbatical Program, by the
MPI (Leipzig, Germany), and by the INFN (Naples, Italy) and PRIN
2002 {\it SINTESI}.

\end{document}